\documentclass[10pt]{article}
\usepackage[utf8]{inputenc}
\usepackage{blindtext}
\usepackage{geometry}
\usepackage{xcolor, float}
 \usepackage{mathtools}
\usepackage{amsthm}
\usepackage{graphicx,amssymb,amsmath,latexsym}
\usepackage{caption}
\usepackage{subcaption}
\usepackage{tikz-cd}
\usetikzlibrary{knots}
\newcommand\Myperm[2][^n]{\prescript{#1\mkern-2.5mu}{}P_{#2}}

\begin{document}
\title{Classification of virtual links by arc shift move}
\author{Aastha Sahore, Komal Negi, Amrender Singh Gill, Prabhakar Madeti}

\theoremstyle{plain}
\newtheorem{theorem}{Theorem}[section]
\newtheorem{lemma}[theorem]{Lemma}
\newtheorem{prop}[theorem]{Proposition}
\newtheorem{corollary}[theorem]{Corollary}

\theoremstyle{definition}
\newtheorem{definition}[theorem]{Definition}
\newtheorem{remark}[theorem]{Remark}
\newtheorem{Example}[theorem]{Example}

\date{}

\maketitle

\begin{abstract}
In this paper, we establish that the arc shift operation on a $n$-component virtual link diagram acts as an unknotting operation when the virtual link is $n$-homogeneous proper, aiding in the classification of \( n \)-component virtual links up to arc shift equivalence. We explore the connection between the arc shift number and the odd writhe of virtual links which are homogeneous proper. Additionally, we identify sequences of virtual link diagrams \( L_n \) for which the upper bound of the arc shift number is exactly \( n \).



\end{abstract}

\section{Introduction}

Virtual links were first introduced by Kauffman~\cite{LH} in $1999$. 
It is an embedding of circular loops into the thickened, compact and oriented surfaces. 
Various knot invariants can be naturally generalized to virtual knot invariants as well as virtual link invariants. To distinguish between two given virtual knots or links, a number of invariants have been introduced. 
The polynomial invariants include the arrow polynomial, index polynomial, multi-variable polynomial, and polynomial invariants of virtual knots and links. 
The numerical invariants includes n-writhe, affine index etc., generalisation of these as virtual link invariants is done by Naoko Kamada et al.\cite{NS}.

Unknotting operations serve as crucial invariants for studying virtual links, such as virtualization of real crossings, and \( H(n) \)-moves, which are explored in the articles~\cite{OM}, and~\cite{DZM}, respectively. However, certain operations such as \(\varXi\)-move, forbidden move and various other local moves that serve as unknotting operations for virtual knots are not applicable to virtual links. 
The classification of virtual links up to the aforementioned local moves is studied in~\cite{MJSW, BPJE}. 
In papers~\cite{AKM, GKP}, the arc shift move is defined and establish that it works as unknotting operations for virtual knots. 
Additionally, the arc shift number for virtual knots is defined and computed for various classes of virtual knots. 
A natural question arises: Is the arc shift operation an unknotting operation for virtual links? In this paper, we demonstrate that not all virtual links can be trivialized by the arc shift move; for instance, the link in Figure~\ref{rew} cannot be trivialized by this move. 
We apply this local move to virtual links and classify them up to arc shift operations based on the parity of their virtual linking numbers, where parity indicates whether the number is odd or even.

This article is structured as follows:
Section 2 presents essential preliminaries about virtual links required for our research. In Section 3, we establish the arc shift as an unknotting operation for a particular class of virtual links, which we term homogeneous proper links, and define the arc shift number for these links. 
 It is interesting to observe that classification of virtual links upto forbidden moves and $\varXi$-move results in an infinite class of virtual links. In contrast, we classify the \(n\)-component virtual links up to arc shift equivalence, demonstrating that they form a finite class. Furthermore, we present sequences of homogeneous proper link diagrams $D_n$, where the upper bound for arc shift number of $D_n$ is $n$.
The article concludes by outlining open problems and potential directions for future research related to these topics.


\section{Preliminaries}
In this article, we consider virtual links as ordered virtual links.
\begin{definition}
An $n$ component virtual link $L$ is the image of an embedding $$f: S_1\cup S_2 \cup \cdots\cup S_n \rightarrow S \times I,$$ where $S_i=S^1$ for all $1\leq i\leq n$ , $S$ is an orientable surface and $I$ is the unit interval.
\end{definition}

\subsection{Virtual link diagrams}
Virtual link diagrams~\cite{LH} are defined as marked generic planar curves, where the markings identify the usual classical crossings, virtual crossings. 

A virtual link is an equivalence class of virtual link diagrams under classical Reidemeister moves(RI, RII, RIII) and virtual Reidemeister moves(VRI, VRII, VRIII,VRIV), together they are known as generalized Reidemeister moves.

\begin{figure}[h]
	\centering
\includegraphics[width=10cm,height=5cm]{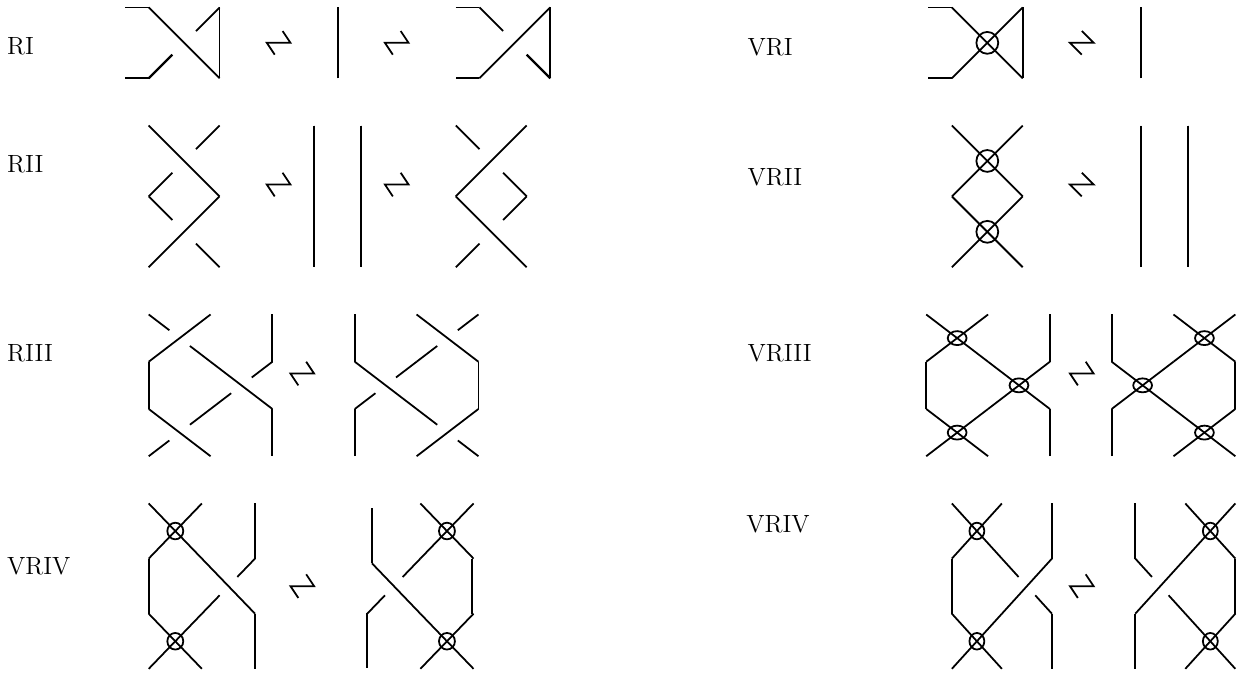}	
  \caption{Generalized Reidemeister moves for virtual link diagrams.}
 \label{rt}
   \end{figure}

\cite{CKS}~It is known that two diagrams represent ambient isotopic virtual link diagrams if and only if they are transformed into
each other by a ﬁnite sequence of generalized Reidemeister moves, as shown in Figure~\ref{rt}. 

\begin{figure}[h]
	\centering
\includegraphics[width=6cm,height=1.5cm]{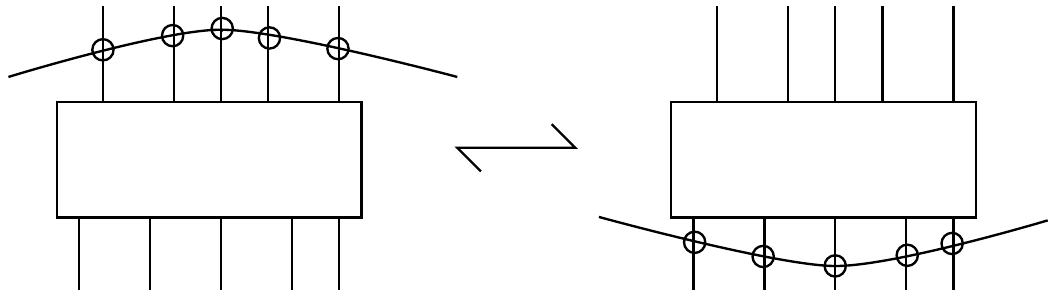}	
  \caption{Detour move.}
 \label{dm}
   \end{figure}

Due to generalized Reidemeister moves, a segment within the virtual link diagram comprises virtual crossings that can be freely move within the plane. When repositioning this segment, the endpoints remain stationary, and any new locations where it intersects the diagram are marked as virtual crossings. This move is referred to as a Detour move, illustrated in Figure~\ref{dm}.

Signs of classical crossings are shown in Figure~\ref{rvs}.
\begin{figure}[h]
    \centering
   \begin{tikzpicture}[scale=1]
	\begin{knot}[xshift=0.5cm]
			\strand[-stealth](1,8) -- (2,9);
			\strand[-stealth](2,8) -- (1,9);
		\end{knot}
			\begin{knot}[xshift=3cm]
			\strand [-stealth](1,8) -- (2,9);
			\strand [-stealth](2,8) -- (1,9);
			\flipcrossings{1};
            \node at (-1.3,8.5) {\text{$c$}};
		  \node at (1.2,8.5) {\text{$c'$}};
			
		\end{knot}
		\end{tikzpicture}
    \caption{Signs of crossing points: $sgn(c)=+1$ and $sgn(c')=-1$.}
    \label{rvs}
\end{figure}
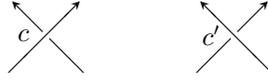
\begin{definition}
Let $L$ be an oriented, $n$-component virtual link where the components are assigned an order: $C_1, C_2, C_3,\ldots, C_n$. Let $C_{j}^{i}$ is the set of crossings of the component $C_{i}$ that overpasses the component $C_{j}$ in $L$. The \textbf{virtual linking number} of $C_i$ over $C_j$ denoted by $L_{j}^{i}$, is defined as\\
$$L_{j}^{i} = \sum_{c\in C_{j}^{i}}{}sgn(c),$$
where $sgn(c)$ is the sign of the crossing $c$.
For 2-component virtual links there are two distinct virtual linking numbers $L_{1}^{2}$ and $L_{2}^{1}$.
\end{definition}

\begin{definition}
    The linking number of oriented, $n$-component virtual link $L$ is defined as follows:
    $$L_{i,j}=\frac{L_{i}^{j} + L_{j}^{i}}{2}$$
\end{definition}
An alternative approach to analyze virtual link diagrams involves using Gauss diagrams.
\subsection{Gauss diagrams}
\begin{figure}[h]
	\centering
\includegraphics[width=16cm]{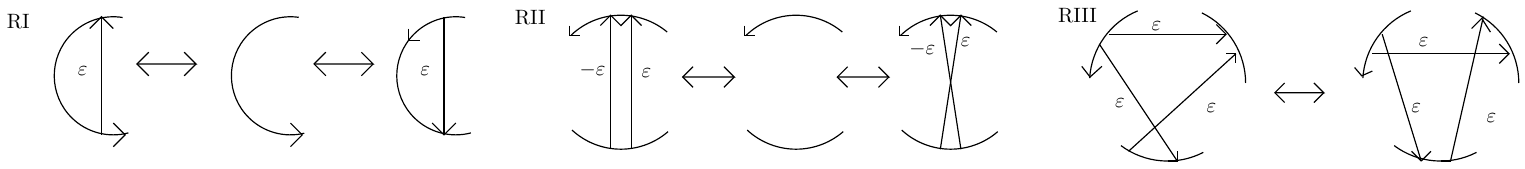}	
\caption{Gauss diagram under generalized Reidemeister moves.}
\label{123}
   \end{figure}

\begin{definition}
\emph{Gauss diagram} $G(D)$ of an $n$-component virtual link diagram $D$ consists of $n$ oriented circles with over(under) passing information in
crossings be presented by directed chords and segments. For a given crossing
$c\in D$ the chord (or segment) in $G(D)$ is directed from over crossing $c$ to
under crossing $c$.
\end{definition}
Two Gauss diagrams correspond to same virtual link if they are equivalent by analogues of Reidemeister moves shown in Figure~\ref{123}.
\begin{Example}
Consider the 2-component virtual link and its corresponding Gauss diagram as shown in Figure~\ref{rew}(a), Figure~\ref{rew}(b), respectively. 
 \begin{figure}[h]
	\centering
\includegraphics[width=8cm]{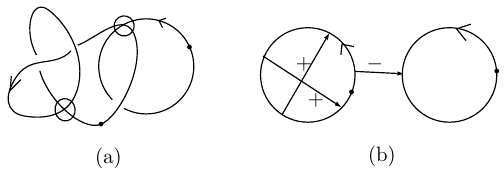}	
  \caption{Gauss diagram for virtual link.}
 \label{rew}
   \end{figure}
\end{Example}

  \subsubsection{Forbidden Moves}
Forbidden moves are not unknotting operation for virtual links and classical links.
Allowing both of the forbidden moves $F_o$ and $F_u$ gives rise to fused isotopy~\cite{LH}. Specifically, two virtual links $L_1$ and $L_2$ are called fused isotopic if $L_2$ can be obtained from $L_1$ by a finite sequence of the $F_o$, $F_u$ moves, and generalized Reidemeister moves for virtual links. Also, known as F-equivalent.

 \begin{figure}[h]
	\centering
\includegraphics[width=12cm,height=2cm]{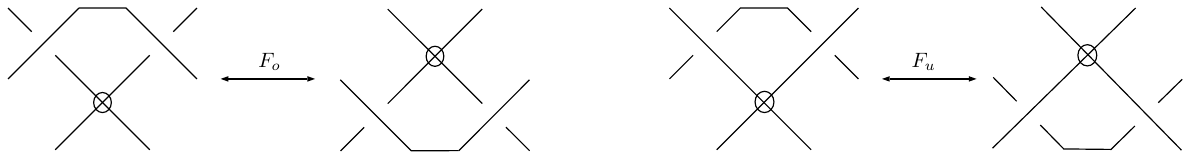}	
  \caption{Forbidden moves for virtual link.}
 \label{fm}
   \end{figure}
   \begin{theorem}~\cite{AE}~\label{ccl}
   A classical link L with n-components is completely determined by the linking numbers of each pair of components under fused isotopy.
   \end{theorem}
     \begin{theorem}~\cite{ABBW}~\label{cwl}
         Two welded links $L$ and $L'$ are F–equivalent if and only if, for all $i , j \in \{1,\ldots, n\}$,
$L_{i}^j = L'^{j}_i$. In other words, fused links are classified by the virtual linking numbers.
     \end{theorem}
\subsubsection{Odd Writhe}
A \textit{2-component even (odd) virtual link} is defined as a virtual link with two components that contains an even (odd) number of mixed crossings.
The definition of odd writhe for $2$-component even virtual links is stated in paper~\cite{MJSW}. 
\begin{definition}\label{odddef}
Consider a $2$-component even virtual link $L=K_1 \cup K_2$ with the virtual link diagram $D=D_1 \cup D_2$.
A self-crossing \( c \) in \( K_i \) is said to be an odd crossing if, following the orientation, we encounter an odd number of real crossings from \( c \) to \( c \). (i.e., odd w.r.t. each component.) The set of all odd crossings of the $i^{th}$ component of $L$ is denoted as $Odd(D_i, D)$.
\end{definition}

\begin{definition}
The sum of signs of all the odd crossings of $i^{th}$-component in a $2$-component even virtual link diagram $D$ is called \emph{odd writhe} of $i^{th}$-component of $D$, denoted by $J(K_i,L)$.\\
$$J(K_1,L) = \sum_{c\in Odd(D_1,D)} sgn(c) \text{, ~~~~~~ } J(K_2,L) = \sum_{c\in Odd(D_2,D)} sgn(c).$$
It is an invariant for $2$- component even virtual links.
\end{definition}

\subsubsection{$\varXi$-move on Gauss diagrams}
\begin{figure}[h]
	\centering
\includegraphics[width=4.5cm,height=2cm]{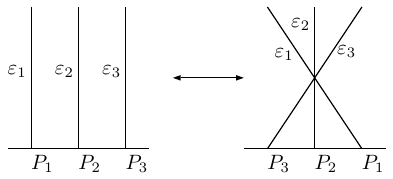}	
  \caption{$\varXi$-move}
 \label{m}
   \end{figure}

Let \( P_1 \), \( P_2 \), and \( P_3 \) be three consecutive endpoints of chords on the same circle in a Gauss diagram. A $\varXi$-move is a deformation that swaps the positions of \( P_1 \) and \( P_3 \) while preserving the signs and orientations of the chords.
The $\varXi$-move on the Gauss diagram is shown in Figure~\ref{m}.
   
If $P_1$ and $P_3$ are connected by a self-chord, then the chord is called a shell.
\begin{definition}~\cite{MJSW}
The two Gauss diagrams $G$ and $G'$ are $\varXi$-equivalent if they are related by a finite sequence of Reidemeister moves RI-RIII and $\varXi$-moves. Two virtual links are $\varXi$-equivalent if their corresponding Gauss diagrams are $\varXi$-equivalent.
\end{definition}

\begin{theorem}~\cite{MJSW}~\label{cvl}
Let $L = K_1 \cup K_2$ and $L'= K'_1 \cup K'_2$
be 2-component odd virtual links. Then the following are equivalent.
\begin{itemize}
    \item[(i)] $L$ and $L'$ are related by a finite sequence of $\varXi$-moves.
\item[(ii)]$L_1^2    = L'^2_1   $ and $L_2^1    = L'^1_2   $.
\end{itemize} 
\end{theorem}
\begin{theorem}~\cite{MJSW}
Let $L = K_1 \cup K_2$ and $L'= K'_1 \cup K'_2$
be 2-component even virtual links. Then the following are equivalent.
\begin{itemize}
    \item[(i)] $L$ and $L'$ are related by a finite sequence of $\varXi$-moves.
\item[(ii)]$J(K_1,L)=J(K_1',L')$, $J(K_2,L)=J(K_2',L')$, and $\overline{F(L)}=\overline{F(L')}$, where $\overline{F(L)}$ is reduced linking class of $L$ which is a refinement of the linking numbers $L^1_2$ and $L_1^2$.
\end{itemize} 
\end{theorem}

\begin{definition}~\cite{MJSW}
    The 2-component virtual links, $L(a_1, a_2, b_1, b_2; k, l)$ and $ M(a_1, a_2, b_1, b_2; k, l)$ for some $a_1, a_2, b_1, b_2, k, l \in \mathbb{Z}$, are defined as virtual links whose Gauss diagrams are represented by $G(a_1, a_2, b_1, b_2; k, l)$ and $ H(a_1, a_2, b_1, b_2; k, l)$, respectively, as shown in Figure~\ref{gh}. 
    A shell-pair is a pair of shells whose four endpoints are consecutive on the same circle. 
    The element $e$ shown in Figure~\ref{gh} is the map $e : \mathbb{Z} \to \{-1, 0, 1\} $ given by 
    $$e(n) = \begin{cases}
~1 , &  n > 0, \\
-1 , & n < 0, \\
~0 , & n = 0.
    \end{cases}$$
    
   \begin{figure}[h]
	\centering
\includegraphics[width=16cm,height=6cm]{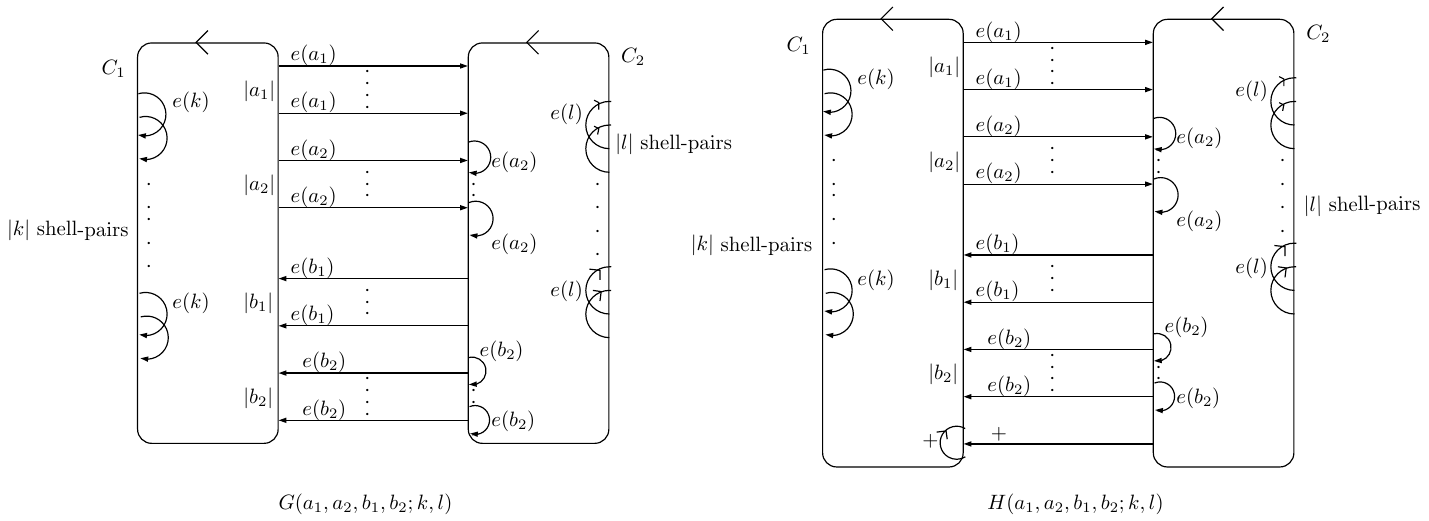}	
  \caption{Gauss diagrams $G(a_1, a_2, b_1, b_2; k, l)$ and $ H(a_1, a_2, b_1, b_2; k, l)$.}
 \label{gh}
   \end{figure}
   \end{definition}
   \begin{prop}~\cite{MJSW}
Any 2-component virtual link $L$ is $\varXi$-equivalent to either
$L(a_1, a_2, b_1, b_2; k, l)$ or $ M(a_1, a_2, b_1, b_2; k, l)$
for some $a_1, a_2, b_1, b_2, k, l \in \mathbb{Z}$.
\end{prop}
\subsubsection{Arc shift move}

\begin{figure}[h]
	\centering
\includegraphics[width=12cm,height=5cm]{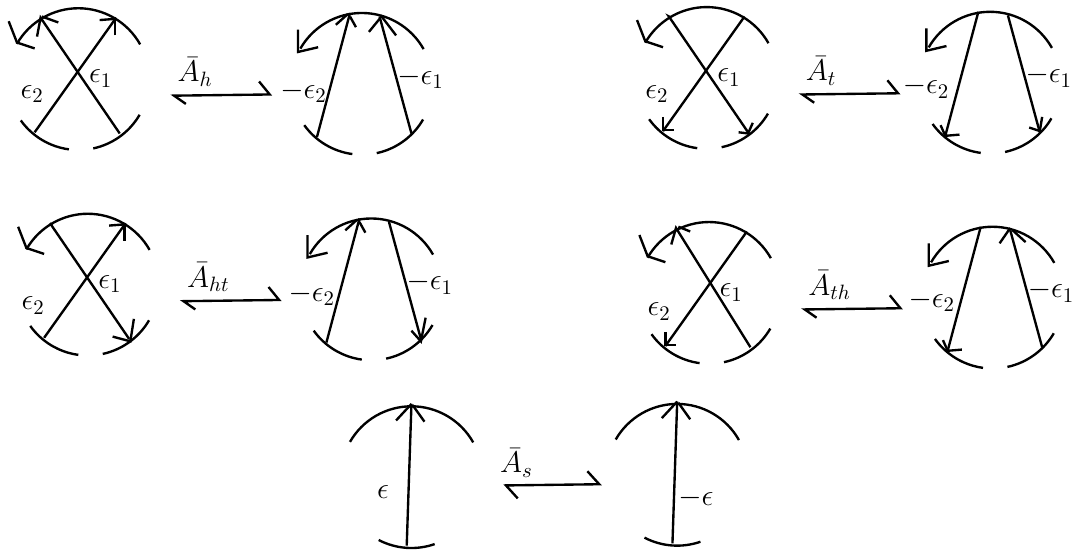}	
  \caption{Gauss diagrams analogues to arc shift moves.}
 \label{g2}
   \end{figure}

 \begin{definition}~\cite{GKP}
Given a virtual knot diagram D, consider an arc $(a,b)$ passing through the pair of crossings $(c_1, c_2)$, where atleast one of the crossings is classical. Then, the arc shift move on the arc $(a,b)$ is shown in Figure~\ref{g2}.
 \end{definition}
   \begin{figure}[h]
	\centering
\includegraphics[width=14cm,height=6cm]{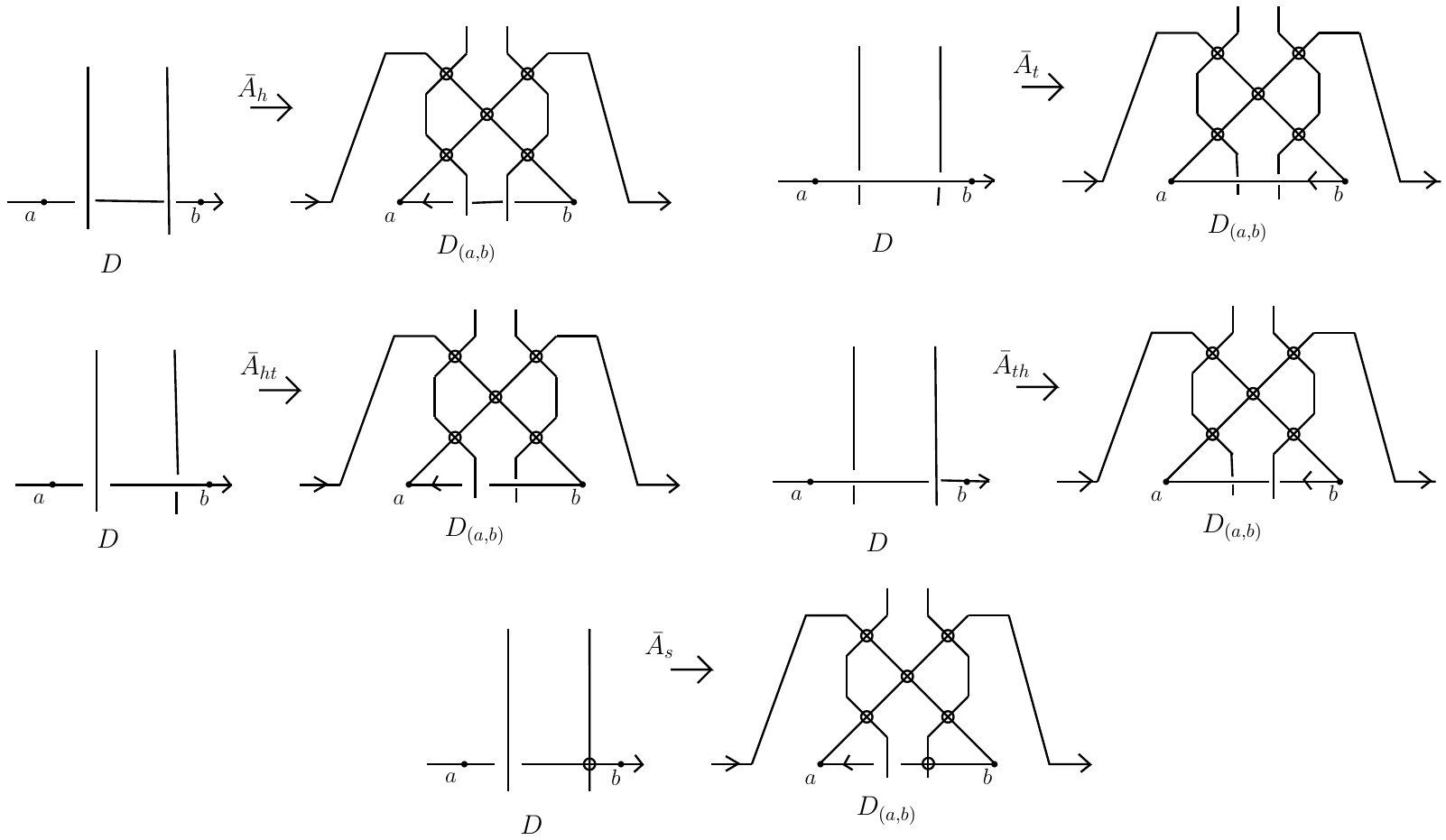}	
  \caption{Arc shift moves on $(a,b)$.}
 \label{asm2}
   \end{figure}
    Moves corresponding to arc shift moves in the Gauss diagram are shown in Figure~\ref{g2}.
    
\begin{theorem}~\cite{GKP}
Every virtual knot diagram $D$ can be transformed into a trivial knot diagram using arc shift moves and generalized Reidemester moves.
\end{theorem}

\begin{definition}
For any virtual knot \( K \), the \textit{arc shift number} is the minimum value of the set \( \{ A(D) \mid D \text{ represents } K \} \), where \( A(D) \) denotes the minimum number of arc shift moves needed to convert a diagram \( D \) into a trivial knot. This is denoted by \( A(K) \).
\end{definition}

\begin{theorem}~\cite{GKP}\label{lowerbound}
Given a virtual knot $K$, arc shift number $A(K)$ and the odd writhe $J(K)$ satisfy the inequality:
$$A(K)\geq |J(K)|/2.$$
Here, $J(K)$ is the sum of signs of the odd crossings of $K$.
\end{theorem}

\begin{theorem}~\cite{AKM}
    For any positive integer $n$, there exist infinitely many virtual knots whose arc shift number is $n$.
\end{theorem}

We have previously noted that arc shift is not an unknotting operation for virtual links. For example, the virtual Hopf link shown in Figure~\ref{rte} cannot be trivialized by an arc shift move. In contrast, the example in Figure~\ref{rst} can be trivialized by an arc shift move. We identify a particular class of virtual links for which this operation serves as an unknotting move, which is discussed in the next section.

\section{Arc shift equivalence for virtual links}

\begin{definition}
 An $n$-component virtual link $L = K_1\cup \cdots \cup K_n$ is said to be $n$-homogeneous proper if $L^i_j   $ and $L^j_i   $ are even for all $1 \leq i\neq j\leq n$.
\end{definition}

 We modify the Definition~\ref{odddef} for $n$-homogeneous proper links. 
\begin{definition}
Consider a $n$-homogeneous proper link $L=K_1 \cup \cdots \cup K_n$ with the diagram $D=D_1 \cup \cdots \cup D_n$.
A self-crossing $c$ in $K_i$ is considered an odd crossing if, following the orientation, we encounter odd number of real crossings from \( c \) to \( c \).(i.e, odd w.r.t. each component.) The set of all odd crossings of the $i^{th}$ component of $L$ is denoted as $Odd(D_i, D)$.
\end{definition}

\begin{definition}
The sum of signs of all the odd crossings of $i^{th}$-component in a homogeneous proper link diagram $D$ is called \emph{odd writhe} of $i^{th}$-component of $D$, denoted by $J(K_i,L)$.\\
$$J(K_i,L) = \sum_{c\in Odd(D_i,D)} sgn(c).$$
It is an invariant for $n$-homogeneous proper links.
\end{definition}
\begin{definition}
    The odd writhe of $n$-homogeneous proper link is $J(L)=\sum_{i=1}^n J(K_i,L)$.
\end{definition}

The arc shift operation for virtual links is the same as the arc shift operation for virtual knots, where crossings could be self or mixed.

\begin{definition}
Two virtual link diagrams $D$ and $D'$ are said to be arc shift-equivalent, denote it by $D \overset{a}{\sim} D'$ if they are related by a finite sequence of generalized Reidemeister moves and arc shift-moves.
\end{definition}

\begin{lemma}~\label{prty}
The parity of virtual linking numbers is invariant under arc shift operation.
\end{lemma}
\begin{proof}
    Let $L = K_1\cup \cdots \cup K_n$ be $n$-component virtual link. 
    Let the arc shift moves be any one of $\overline{A}_{h}$, $\overline{A}_{t}$, $\overline{A}_{ht}$, ${\overline{A}}_{th}$, $\overline{A}_s$  that we apply on the pair of crossings say $(c_1,c_2)$.
    Let $L' = K'_1\cup \cdots \cup K'_n$ be the virtual link obtained by applying one of the arc shift moves mentioned above. It is sufficient to prove that parity of $L_i^j $ equals to the parity of $L'^j_i  $.
    Here the following cases arises:
    \begin{itemize}
        \item[] Case(I) When both $c_1$ and $c_2$ does not belongs to $C_i^j$. There will be no effect on virtual linking numbers.
        \item[]Case(II) When either $c_1$ or $c_2$ belongs to $C_i^j$, but not both. W.L.O.G., let $c_1\in C_i^j$ with sign of $c_1$ is $\varepsilon$. After applying arc shift move, sign of $c_1$ changes to $-\varepsilon$ and hence $$L'^j_i  =L_i^j -2\varepsilon.$$
        Therefore the parity of $L_i^j $ is same as the parity of $L'^j_i  $. 
\item[]Case(III) When both $c_1$ and $c_2$ are in $C_i^j$. Let sign of $c_1$ be $\varepsilon_1$ and sign of $c_2$ be $\varepsilon_2$. After applying arc shift move, sign of $c_1$ and $c_2$ changes to $-\varepsilon_1$, and $-\varepsilon_2$, respectively. Hence, $$L'^j_i  =L_i^j -2(\varepsilon_1+\varepsilon_2).$$
Therefore the parity of $L_i^j $ is same as the parity of $L'^j_i  $.
    \end{itemize}
\end{proof}
\begin{corollary}~\label{cp}
    If one of the parity of virtual linking number of virtual link $L$ is odd then $L$ cannot be unlinked by arc shift operation.
\end{corollary}
\begin{Example}
Consider a virtual link shown in Figure~\ref{rte}. The virtual linking numbers of this virtual link are $L_1^2   =1$ and $L_2^1   =0$. Since the parity of one of the virtual linking numbers is odd, the link cannot be trivialized by the arc shift operation. 
 \begin{figure}[h]
	\centering
\includegraphics[width=4cm,height=1.5cm]{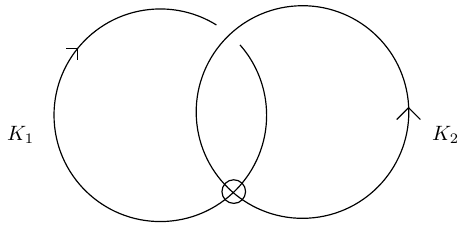}	
  \caption{Virtual Hopf link}
 \label{rte}
   \end{figure}
   \end{Example}
   
\begin{figure}[ht]
	\centering
\includegraphics[width=8cm]{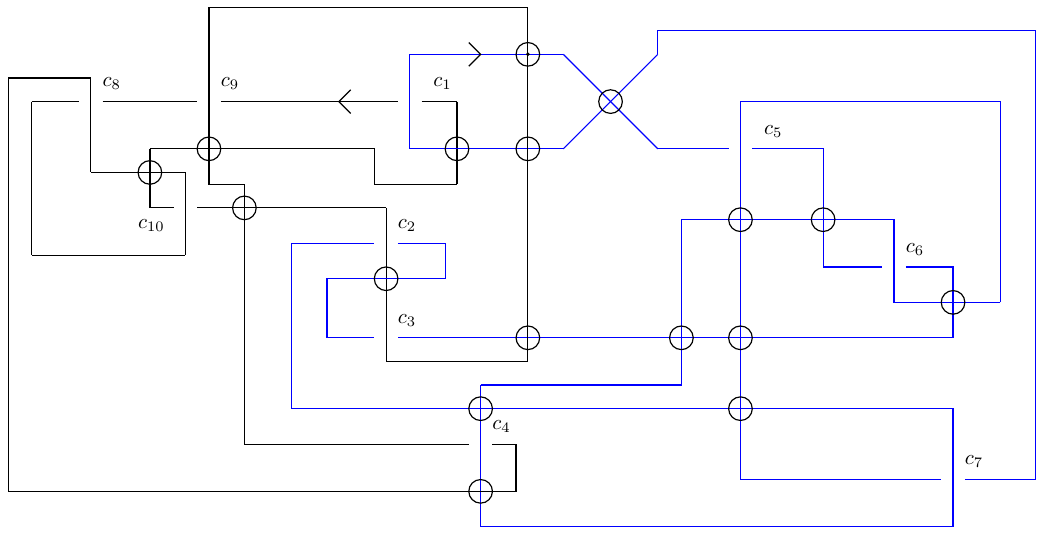}	
  \caption{2-homogeneous proper link.}
 \label{gdm}
   \end{figure}

\begin{figure}[!]
	\centering
\includegraphics[width=12cm]{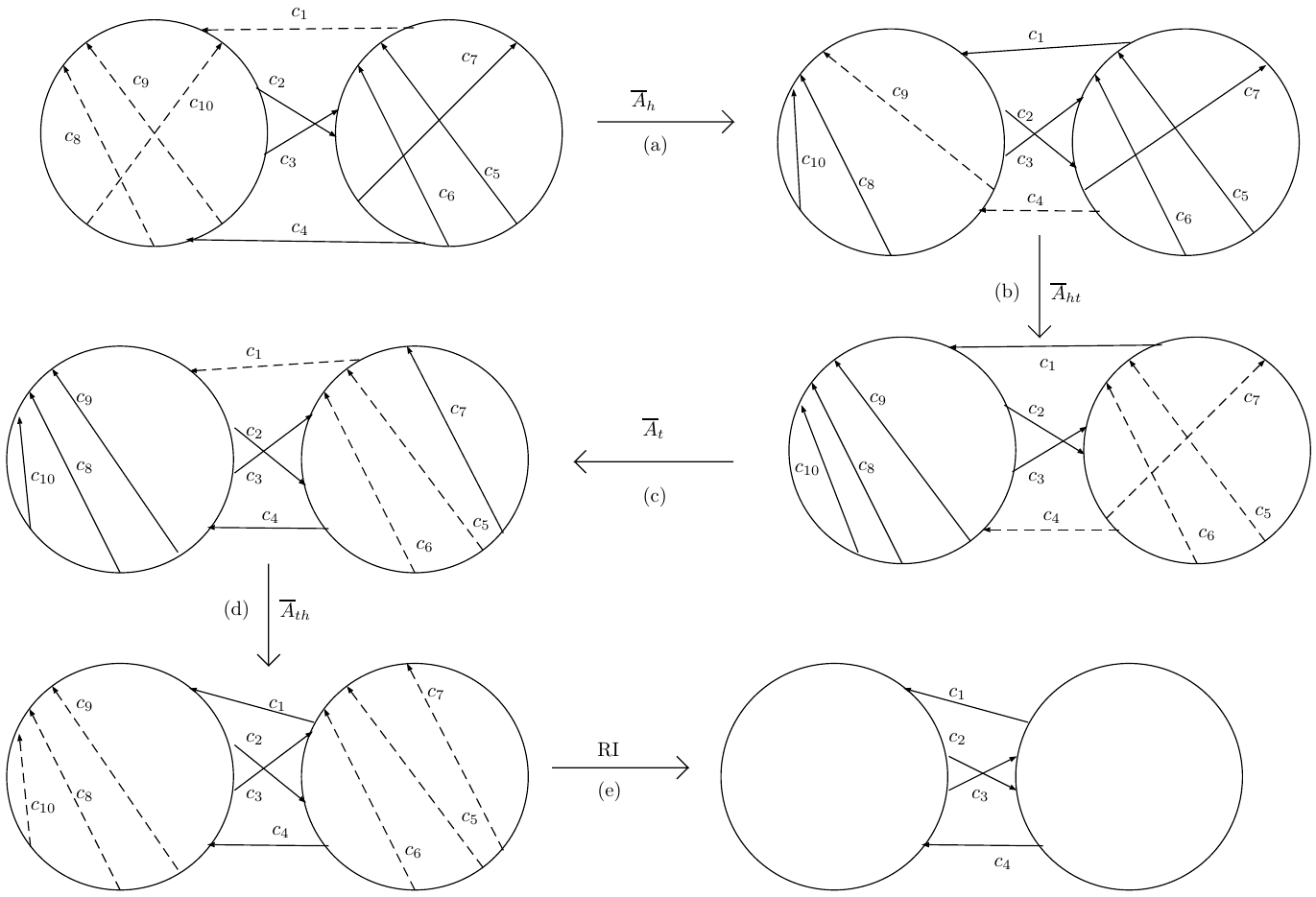}	
  \caption{Parallel alignment of self-chords of a Gauss diagram.}
 \label{gdl}
  \end{figure}

Consider the 2-componenet virtual link shown in Figure~\ref{gdm}.     
We will unlink this example by applying arc shifts on this Gauss diagram, as shown in Figure~\ref{gdl}. We will demonstrate the different stages of this process in the following remarks.

\begin{remark}~\label{rem1}
   For the transformation of self-chords of a Gauss diagram into parallel alignment, we begin by applying the arc shift move \( \overline{A}_h \) to the pairs \((c_1,c_{10})\), then \((c_9,c_{10})\), and finally \((c_8,c_{10})\), as shown in Figure~\ref{gdl}(a). Next, the arc shift move \( \overline{A}_{ht} \) is applied to the pair \((c_9,c_4)\), which aligns the self-chords of the first component in parallel, as seen in Figure~\ref{gdl}(b). 
The next step is to align the self-chords of the second component. This is achieved by applying the arc shift move \( \overline{A}_t \) to the pairs \((c_7,c_4)\), \((c_7,c_6)\), and \((c_7,c_5)\) consecutively. Then, we apply the arc shift move \( \overline{A}_{th} \) to the pairs \((c_1,c_5)\) and then \((c_1,c_6)\), as shown in Figures~\ref{gdl}(c)-(d). At this stage, all self-chords of both components are parallel. Finally, using the RI move, these chords can be removed from the Gauss diagram, as demonstrated in Figure~\ref{gdl}(e). This procedure helps in removing the self chords of the given Gauss diagram.
  
\end{remark}
  
\begin{remark}~\label{rem2}
Consider a Gauss diagram consisting of mixed chords only. Our objective is to align these mixed chords in parallel. This procedure is illustrated in Figure~\ref{gdl2}, where the arc shift move \( \overline{A}_h \) is applied to the chord pair \((c_2, c_3)\) in Figure~\ref{gdl2}(a). The resulting Gauss diagram, after the operation, is depicted in Figure~\ref{gdl2}(b).

\begin{figure}[ht]
	\centering
\includegraphics[width=10cm]{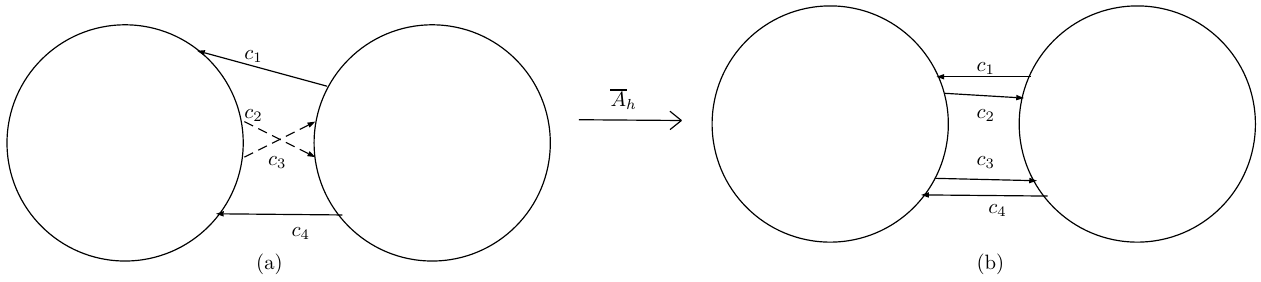}	
  \caption{Parallel alignment of mixed-chords of a Gauss diagram.}
 \label{gdl2}
  \end{figure}
\end{remark}

  \begin{remark}~\label{rem3}
The Gauss diagram consists solely of mixed chords, which are already parallel. Our goal is to arrange the chords in such a way that it facilitates the easy application of the RII move. We apply the arc shift move \( \overline{A}_{th} \) consecutively to the pairs \((c_2,c_1)\) and \((c_3,c_1)\), which arranges all the mixed chords running from left to right above those running from right to left. The final result of this rearrangement is illustrated in Figure~\ref{gdl3}(a).
\begin{figure}[ht]
	\centering
\includegraphics[width=15cm]{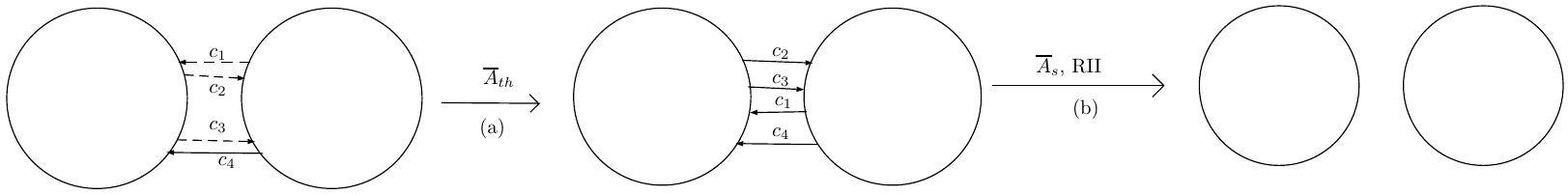}	
  \caption{Trivializing the Gauss diagram.}
 \label{gdl3}
  \end{figure}

In the final step, the arc shift move \( \overline{A}_s \) is applied to the mixed crossings based on their signs, ensuring that consecutive pairs have opposite signs. Subsequently, by applying the RII move, all the chords can be removed, as the mixed chords form pairs, as depicted in Figure~\ref{gdl3}(b).
\end{remark}

In this subsection, we will prove that the arc shift is an unknotting operation for $2$-homogeneous proper links. Later we will extend this result to $n$-homogeneous proper links using induction.
   
\subsection{Main results}
 \begin{theorem}~\label{main}
Arc shift is an unknotting operation on a diagram of 2-component virtual link $L = K_1\cup K_2$ if and only if $L$ is a 2-homogeneous proper link.
   \end{theorem}
   \begin{proof}
   $(\Rightarrow)$
For the forward direction, we will proceed by contradiction. Suppose, contrary to our assumption, that the given 2-component virtual link \(L\) is not a 2-homogeneous proper link. Then the parity of at least one of the virtual linking numbers must be odd. By Corollary~\ref{cp}, this implies that \(L\) cannot be unknotted by the arc shift operation.

   $(\Leftarrow)$
       Let \( L = K_1 \cup K_2 \) be a 2-homogeneous proper link. Then \( L^1_2 \) and \( L^2_1 \) are both even. Consequently, the total number of mixed crossings is even.

     Let \( D = D_1 \cup D_2 \) denote the link diagram corresponding to \( L \), and let \( G(D) = G(D_1 \cup D_2) \) be the Gauss diagram associated with \( D \). Trivializing \( G(D) \) using arc shift operations is equivalent to trivializing \( D \) through arc shift moves.

Initially, we align all self-chords parallel within the same component to resolve self-crossings. This alignment is achieved by arc shift operations. An arc shift operation moves the head of a selected self-chord along with an adjacent chord, while simultaneously switching the signs of both chords. The adjacent chord may be either a self-chord or a mixed chord, and this operation applies to both types. The process is repeated for each self-chord encountered, following the orientation, until a Gauss diagram is obtained where no chord lies between the head and tail of the initially chosen self-chord. Repeating this process for all self-chords yields a Gauss diagram where no two chords intersect, resulting in the required parallel self-chords diagram.

Since all the chords belonging in the same component can be removed by RI move. As a result, we get a Gauss diagram only comprising chords representing mixed crossings. This process is depicted in Remark~\ref{rem1}.

Next, we align all mixed chords parallel using arc shift operations. Begin by selecting the uppermost chord of the first component that does not occupy the same uppermost position in the second component. Suppose this chord is in the $n^{th}$ position on the second component. Perform an arc shift operation between the chords in the $n^{th}$ and $(n-1)^{th}$ positions on the second component. This shifts the $n^{th}$ chord to the $(n-1)^{th}$ position. Continue this process until the selected chord occupies the same position in both components. Repeating this for each chord sequentially will ensure all mixed chords are aligned in parallel. This process is depicted in Remark~\ref{rem2}, through an example.

This Gauss diagram can be further rearranged using arc shift moves so that all mixed chords running from left to right are positioned above those running from right to left. To achieve this, begin by selecting the first mixed chord that runs from left to right, assuming it is positioned as the $n^{th}$ mixed chord from the top. Apply the arc shift operation twice to the $n^{th}$ and $(n-1)^{th}$ chords, which results in switching their positions. By repeating this process $(n-1)$ times, the selected chord moves to the topmost position among all the mixed chords. Then, repeat the same process for the next parallel chord from left to right, using arc shift operations to move it to the second-highest position. Continue this process until all chords running from left to right are positioned above those running from right to left.

   Note that an arc shift move can change the sign of chords but cannot invert them; that is, the count of heads and tails on a component remains unchanged after applying an arc shift operation. Given that \( L \) is a 2-homogeneous proper link, the number of mixed chords is even. Consequently, we can pair mixed chords in such a way that adjacent chords have opposite signs, this can be achieved using the move \( \overline{A}_s \). Subsequently, by sequentially applying the RII move on the Gauss diagram to each paired chord with opposite signs, we can trivialize \( G(D) \). This process is depicted in Remark~\ref{rem3}.

   \end{proof}
   
\begin{Example}~\label{e1}
Consider a $2$-homogeneous proper link is shown in Figure~\ref{ctrefoil}. The virtual linking numbers are $L^2_1=2$ and $L^1_2=0$.
 \begin{figure}[ht]
	\centering
\includegraphics[width=10cm]{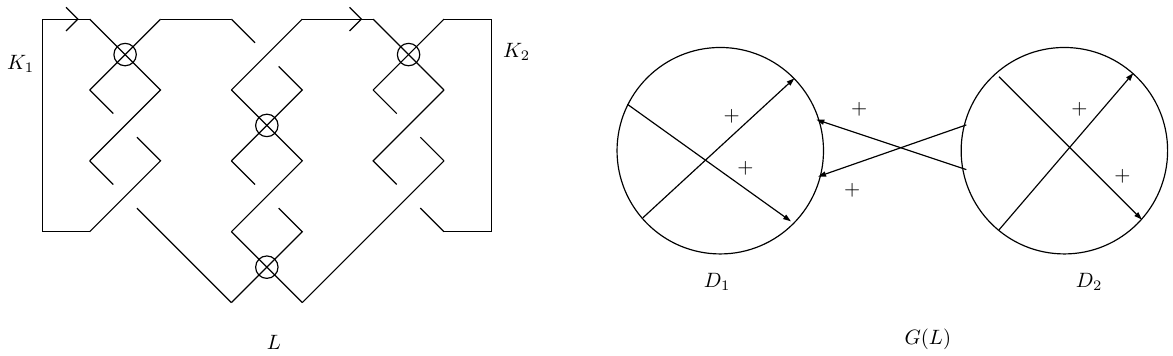}	
  \caption{The $2$-homogeneous proper link diagram $L$.}
 \label{ctrefoil}
   \end{figure}
   
We present an algorithm to unknot the given 2-component virtual link using arc shift moves. First, we apply the arc shift operation on self-crossings, specifically \(\overline{A}_{ht}\) and \(\overline{A}_{th}\) in this case, as shown in Figure~\ref{sf}(a). By applying the RI move, we can eliminate the self-chords as illustrated in Figure~\ref{sf}(b). We are then left with only mixed crossings. By utilizing the \(\overline{A}_{s}\) move, we can change the sign of one chord, as depicted in Figure~\ref{sf}(c). Subsequently, we can apply the RII move to trivialize the diagram.

\begin{figure}[ht]
	\centering
\includegraphics[width=10cm]{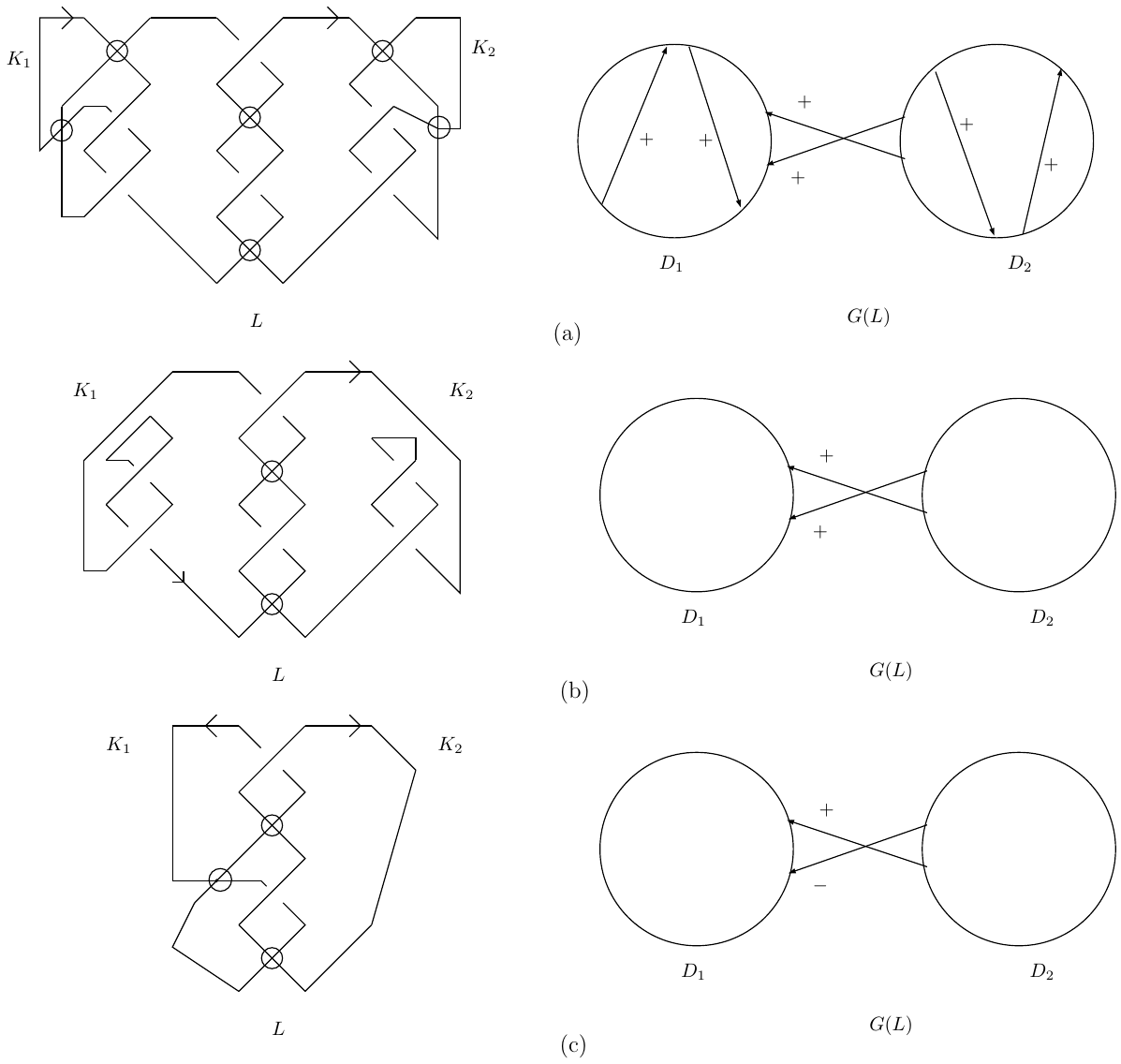}	
  \caption{Unknoting sequence of $2$-homogeneous proper link diagram $L$.}
 \label{sf}
   \end{figure}
\end{Example}

   \begin{corollary}
Arc shift is an unknotting operation on a diagram of the virtual link $n$ of components $L = K_1\cup \cdots \cup K_n$ if and
only if $L$ is a $n$-homogeneous proper link.
 \end{corollary}
 \begin{proof}
    For the proof, we will apply induction on the number of components of a virtual link \( L \). For \( n=1 \), the virtual link is a virtual knot, and for virtual knots, the arc shift is an unknotting operation. So, we begin with \( n=2 \), and by the Theorem~\ref{main}, this case holds true. Now, assume that for \( n=k \), any virtual link with \( k \) components can be unlinked, provided that \( L^i_j \) and \( L^j_i \) are even for \( 1 \leq i \neq j \leq k \). We need to prove this for \( n=k+1 \).

First, select the two leftmost components of the virtual link. Apply the arc shift on the first component to make its self-chords parallel, as done in the Theorem~\ref{main}. By applying the RI move, the self-chords of the first component can be removed. Next, use arc shifts to parallel all the mixed chords between the first two components, ensuring that there are no self-chords of the second component located between the mixed chords of the first two components. Then, apply the arc shift \( \overline{A}_s \) to adjust the signs of the mixed crossings as needed, and finally, use the RII move to remove all the mixed chords between the first two components. This is possible since both \( L^1_2 \) and \( L^2_1 \) are even.

At this point, we are left with a virtual link having \( k \) components. By inductive assumption, this \( k \)-component virtual link can be unlinked.
 \end{proof}
Now, we will explore the class of other virtual links that are not $n$-homogeneous proper links, up to the arc shift move.
 \begin{definition}
 A special class of 2-component virtual link, denoted by $L(p,q)$ is defined as $p, q \in \mathbb{Z}$ must satisfy the following conditions:
      \begin{itemize}
          \item[(i)] Cardinality of the set $C_1^2(L)$ is $|p| $ , and $sgn(c)=e(p)$ for each $c \in C_1^2(L)$, 
          \item[(ii)] Cardinality of the set $C_2^1(L)$ is $|q|$, and $sgn(c)=e(q)$ for each $c \in C_2^1(L)$.
      \end{itemize}
 \end{definition}
\begin{prop}\label{c2cl}
     Let $L = K_1 \cup K_2$ be a $2$-component virtual link. Then $L$ is arc shift-equivalent to $L(m, n)$, where $m=p \text{ mod } 2$ and $n=  q \text{ mod }2$ such that $p = L_1^2   $ and $q = L_2^1   $. 
\end{prop}
\begin{proof}
  If $L$ is 2-homogeneous proper then $L^1_2   $ and $L^2_1   $ are even. By Theorem~\ref{main}, $L\overset{a}{\sim} L(0,0)$. Suppose $L$ is not 2-homogeneous proper, meaning either $L^1_2   $ is odd or $L^2_1   $ is odd. This leads to two cases:

Case(I) If both \( L^1_2    \) and \( L^2_1    \) are odd, we will ultimately obtain the Gauss diagram of \( L(1,1) \) as shown in Figure~\ref{2cvl}(c), even after applying the arc shift moves and RII move.

Case(II) If either \( L^1_2    \) or \( L^2_1    \) is odd, but not both, the number of mixed crossings will be odd. Consequently, the Gauss diagram can be transformed into either \( L(1,0) \) or \( L(0,1) \), using the arc shift moves and RII move,  as depicted in Figure~\ref{2cvl}(a), and ~\ref{2cvl}(b), respectively.

\begin{figure}[h]
	\centering
\includegraphics[width=13cm,height=2cm]{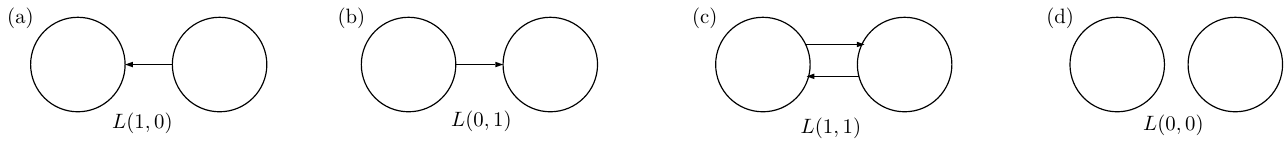}	
  \caption{Gauss diagrams of $2$-component virtual links upto arc shift move.}
 \label{2cvl}
   \end{figure}

\end{proof}
\begin{remark}
From Proposition~\ref{c2cl}, it is evident that every 2-component virtual link is equivalent to one of the virtual links, namely $L(0,0)$, $L(1,0)$, $L(0,1)$, and $L(1,1)$ upto arc shift equivalence.

Since each two of $L(0,0)$, $L(1,0)$, $L(0,1)$, and $L(1,1)$ have atleast one of the virtual linking number's parity is different from the other, hence from Lemma~\ref{prty}, each of them are mutually inequivalent upto arc shift move.

 \begin{figure}[h]
	\centering
\includegraphics[width=12cm,height=2cm]{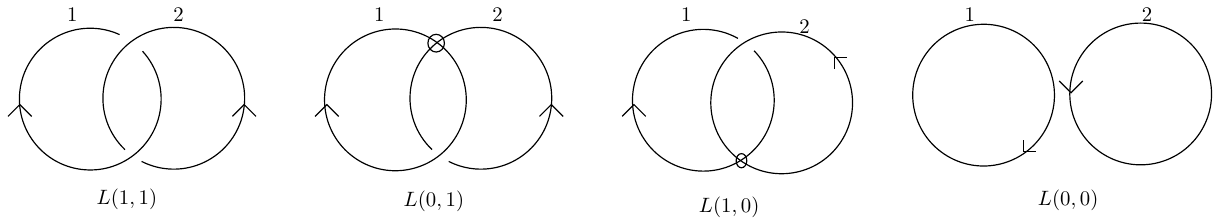}	
  \caption{$2$-component virtual links upto arc shift move.}
 \label{componets}
   \end{figure}
   As a consequence it becomes clear that there are exactly four 2-component virtual links, $L(0,0)$, $L(1,0)$, $L(0,1)$, and $L(1,1)$ upto arc shift equivalence. These four virtual links and their corresponding Gauss diagram are shown in Figure~\ref{componets} and Figure~\ref{2cvl}, respectively. 
\end{remark}

\begin{Example}
       For each positive integer $n$, we can construct a class of 2-component virtual link denoted as $L(2n-1)$ such that the given link has $2n-1$ classical crossings as shown in Figure~\ref{class1}. The virtual linking numbers of $L(2n-1)$ are $L_1^2=0$ and $L_2^1=2n-1$. Every member of the class is arc shift equivalent to the virtual link $L(1,0)$.
\begin{figure}[h]
	\centering
\includegraphics[width=10cm,height=4cm]{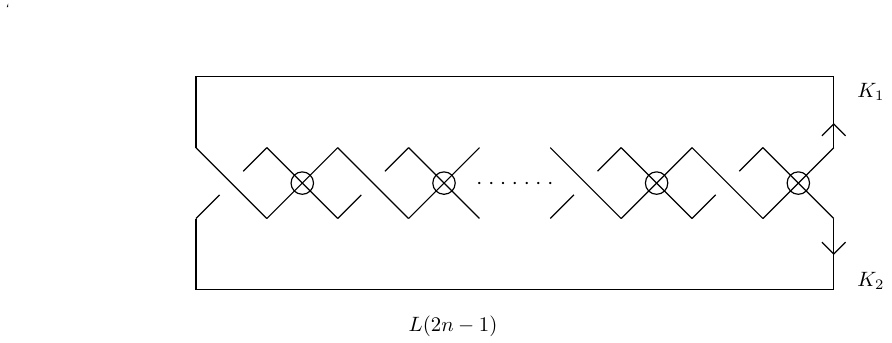}	
  \caption{A class of virtual link $L(2n-1)$.}
 \label{class1}
   \end{figure}
\end{Example}
\begin{definition}
    Let $D$ be a virtual knot diagram and $D^*$ is said to be mirror image of $D$ if $D^*$ is obtained from $D$ by changing the crossing type at every classic crossing point.
\end{definition}
\begin{prop}
Let $L = K_1 \cup K_2$ be a $2$-component virtual link diagram, and $L^* = K^*_1 \cup K^*_2$ be the mirror image of $L$. Then $L \overset{a}{\sim} L^*$ if and only if the parity of $L_1^2   $ and $L_2^1   $ are the same. 
\end{prop}
\begin{proof} 
Let $L^{*}$ be a mirror image of $L$, and $x^*$ is a real crossing point of $L^*$ corresponding to $x$ in  $L$ such that $sgn(x)=-sgn(x^*)$.
If $x\in C_{2}^{1}(L)$, then $x^* \in C_{1}^{2}(L^*)$.
Hence $L_{2}^{1}   =-L_{1}^{*2}$.
Similarly, if $x \in C_{1}^{2}(L)$ then $x^* \in C_{2}^{1}(L^*)$.
Hence $L_{1}^{2}   =-L_{2}^{*1}$. Therefore, by Corollary~\ref{c2cl}, \begin{equation}\label{1}
     L\overset{a}{\sim}L((m\mod 2), (n\mod 2)) \text{ and } L^*\overset{a}{\sim}L((n\mod 2), (m\mod 2)),
\end{equation}
where  $m = L_1^2   $ and $n = L_2^1   $.

$(\Rightarrow)$
Forward part follows from equation~(\ref{1}) and Lemma~\ref{prty}.

$(\Leftarrow)$
For the converse part, let the parity of $L_1^2   $ and $L_2^1   $ are the same. From equation~(\ref{1}), $m \equiv n \mod 2$. The following holds,
\begin{equation}\label{2}
     L\overset{a}{\sim}L((n\mod 2), (n\mod 2)) \text{ and } L^*\overset{a}{\sim}L((n\mod 2), (n\mod 2)),
\end{equation}
From equation~(\ref{2}), we have $L \overset{a}{\sim} L^*$.

\end{proof}
   \begin{theorem}
       If $L=K_1 \cup \cdots \cup K_n$ and $L'=K_1' \cup \cdots \cup K'_n$ are two virtual links then the following are equivalent:
       \begin{itemize}
           \item[(i)] $L$ and $L'$ are related by finite sequence of arc shift moves
           \item[(ii)] $L^i_j   =(L'^i_j   \mod 2)$ and $L^j_i   =( L'^j_i    \mod 2)$ for $1 \leq i \neq j\leq n$.
       \end{itemize}
   \end{theorem}
\begin{proof}
$(i) \implies (ii)$ It can be proved by Lemma~\ref{prty}. 
 
$(ii) \implies (i)$ Let $i\neq j$, the parity of $L^i_j   $ and $L^j_i   $ is equal to parity of $L'^i_j   $ and $L'^j_i   $, respectively.
 By the reducing algorithm, depicted in Theorem~\ref{main}, both $L$ and $L'$ can be reduced by arc shift operation to same virtual link $L''$ such that the parity of virtual linking numbers of $L''$ is same as the parity of virtual linking numbers of $L$ and $L'$. We have $L\overset{a}{\sim} L''$ and $L''\overset{a}{\sim} L'$, this implies $L\overset{a}{\sim} L'$.
\end{proof}
      
   \begin{corollary}
     Let $L = K_1 \cup \cdots \cup K_n$ be a $n$-component virtual link diagram, and $L^* = K^*_1 \cup \cdots \cup K^*_n$ be the mirror image of $L$. Then $L \overset{a}{\sim} L^*$ if and only if the parity of all virtual linking numbers of $L$ are the same. 
\end{corollary}
  \begin{remark}
 By Proposition~\ref{c2cl}, we conclude that, up to arc shift equivalence, there are four distinct 2-component virtual links, as illustrated in Figure~\ref{componets}. This result can be generalized for $n$-component virtual links, where there are $2^{(\Myperm{2})}$ distinct $n$-component virtual links, up to arc shift moves.
   \end{remark}
   As discussed in the preliminaries, Theorem~\ref{ccl} addresses classical links up to forbidden moves, Theorem~\ref{cwl} deals with welded links up to forbidden moves, and Theorem~\ref{cvl} examines virtual links up to the \(\varXi\)-move. Each of these cases results in infinite classes of classical links, welded links, and virtual links, respectively. In contrast, the virtual links up to the arc shift move yields a finite class, distinguishing it from the infinite nature of the other cases.

   The following proposition provides a condition under which the class of 2-component virtual links up to the \(\varXi\)-move becomes equivalent to the finite class of 2-component virtual links up to the arc shift move.
   \begin{prop}
    Given 2-component virtual links $L(a_1, a_2, b_1, b_2; k, l)$ and $M(a_1, a_2, b_1, b_2; k, l)$ then the following are true:
    \begin{itemize}
        \item[(i)] If $a_1+a_2\equiv 0 \mod 2$ and $b_1+b_2\equiv 0 \mod 2$, then $L(a_1, a_2, b_1, b_2; k, l) \overset{a}{\sim} L(0,0)$ and $M(a_1, a_2, b_1, b_2; k, l) \overset{a}{\sim} L(1,0)$
        \item[(ii)] If $a_1+a_2\equiv 0 \mod 2$ and $b_1+b_2\equiv 1 \mod 2$, then $L(a_1, a_2, b_1, b_2; k, l) \overset{a}{\sim} L(1,0)$ and $M(a_1, a_2, b_1, b_2; k, l) \overset{a}{\sim} L(0,0)$
        \item[(iii)] If $a_1+a_2\equiv 1 \mod 2$ and $b_1+b_2\equiv 0 \mod 2$, then $L(a_1, a_2, b_1, b_2; k, l) \overset{a}{\sim} L(0,1)$ and $M(a_1, a_2, b_1, b_2; k, l) \overset{a}{\sim} L(1,1)$
        \item[(iv)] If $a_1+a_2\equiv 1 \mod 2$ and $b_1+b_2\equiv 1 \mod 2$, then $L(a_1, a_2, b_1, b_2; k, l) \overset{a}{\sim} L(1,1)$ and $M(a_1, a_2, b_1, b_2; k, l) \overset{a}{\sim} L(0,1)$
    \end{itemize}
\end{prop}
   \subsection{Arc shift number}
   \begin{definition}
 For any given $n$-homogeneous proper link $L$, its \textit{arc shift number} denoted as $A(L)$ is determined as the minimum value of the set $\{ A(D) | D \text{ represents } L \}$, where $A(D)$ signifies the minimum number of arc shift moves required to transform a diagram $D$ into a trivial link.
\end{definition}

\begin{theorem}
Arc shift number is an invariant for $n$-homogeneous proper link $L$. 
\end{theorem}

\begin{Example}
    Consider the virtual link \( L \) depicted in Figure~\ref{rst}. Since $L_{1}^{2}   =0$, and $L_{2}^{1}   =2$ both are even, hence by Lemma~\ref{prty}, it is arc shift equivalent to trivial virtual link. Applying the arc shift operation to the arc \((a,b)\) allows us to unknot the link \( L \) by subsequently performing a detour move followed by Reidemeister moves (RII and VRII). This indicates that the arc shift number for virtual links is at most \( 1 \). Furthermore, since the link \( L \) is non-trivial, at least one arc shift is necessary to unknot the virtual link. Thus, the arc shift number for this virtual link is exactly \( 1 \).
\begin{figure}[h]
	\centering
\includegraphics[width=14cm]{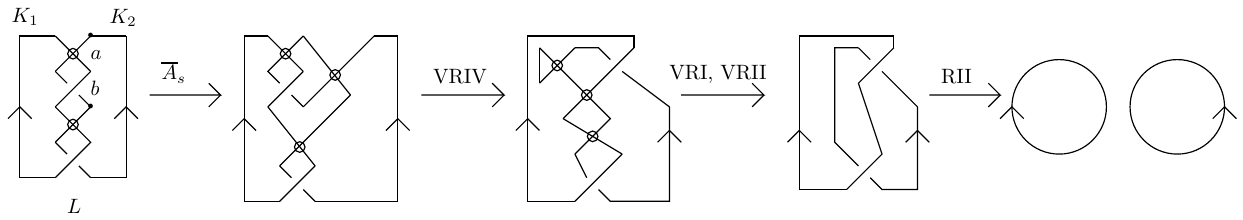}	
  \caption{Arc shift on virtual link}
 \label{rst}
   \end{figure}
\end{Example}

\begin{theorem}~\label{J(K)}
    If $L=K_1 \cup \cdots \cup K_n$ and $L'=K_1' \cup \cdots \cup K'_n$ are $n$-homogeneous proper link diagrams that differ by an arc shift move, then $J(K_i', L') = J(K_i,L)$ or $J(K_i', L') = J(K_i,L)\pm1$ or $J(K_i', L') = J(K_i,L)\pm2$, for all $i$.
\end{theorem}

\begin{proof}
   Let the arc shift moves be any one of $\overline{A}_{h}$, $\overline{A}_{t}$, $\overline{A}_{ht}$, ${\overline{A}}_{th}$, $\overline{A}_s$  that we apply on the pair of crossings say $(c_1,c_2)$. Let $c_1'$, $c_2'$ be the corresponding chords after applying the arc shift move. If both the crossings are self crossings, then the change will take place w.r.t. a single component. Then by Theorem~\ref{lowerbound}, there is a change in sign and the parity of each crossings and the odd writhe will be either $J(K_i, L) = J(K_i',L)$ or $J(K_i, L) = J(K_i',L)\pm2$.
  In case of arc shift on a self crossing and a mixed crossing, the sign of self crossing and the parity will change, and the sign of mixed crossing will change while keeping the parity unchanged, because it does not contribute in $J(K_i,L)$. We will discuss the two cases based on the parity of self crossing $c_1$.\\
Case (I) When the crossings $c_1$ is even. Then, $c_1'$ is odd contribute in $J(K_i',L')$. Therefore, 
  \begin{align*}
  J(K_i',L') &= J(K_i,L)+sgn(c_1')\\
        &= J(K_i,L)+(-\varepsilon_1)\\
        &= J(K_i,L)-\varepsilon_1
  \end{align*}
Case (II) When the crossings $c_1$ is odd. Then, $c_1'$ is even, does not contribute in $J(K_i',L')$. Therefore, 
  \begin{align*}
      J(K_i',L') &= J(K_i,L)-sgn(c_1)\\
            &= J(K_i,L)-\varepsilon_1
  \end{align*}

From the above cases, we have $J(K_i',L') = J(K_i,L)\pm1$.  

Therefore, if $L$ and $L'$ are two virtual link diagrams that differ by an arc shift move, then $J(K_i',L') = J(K_i,L)$ or $J(K_i',L') = J(K_i,L)\pm1$ or $J(K_i',L') = J(K_i,L)\pm2$.
\end{proof}

\begin{theorem}~\label{Lowerbd}
For a $n$-component homogeneous proper link $L$, arc shift number $A(L)\geq J(L)/2.$ 
\end{theorem}

\begin{proof}
    Consider a virtual link $L$. Let the arc shift number of $L$ is $n$. Let $L_0\rightarrow L_1\rightarrow L_2\rightarrow\cdots \rightarrow L_t$ be the sequence of generalized Reidemeister moves and arc shift moves applied to unlink $L.$ $L_t$ is a trivial link and $L_i$ is obtained from $L_{i-1}$ either by arc shift or by generalized Reidemeister move. Exactly $n$ terms in the sequence corresponds to the arc shift moves. We have,
    \begin{align}
        |J(L_t) - J(L_0)| &= |J(L_t)-J(L_{t-1})+
J(L_{t-1})-J(L_{t-2})+....+J(L_1)-J(L_0)|\\
         &\leq |J(L_t)-J(L_{t-1})|+|
J(L_{t-1})-J(L_{t-2})|+....+|J(L_1)-J(L_0)|\label{3}
\end{align}
where, $J(L_t) = 0$ and $J(L_0)=J(L)$. There are exactly $n$ terms in the inequality~(\ref{3}) which correspond to arc shift moves, while the remaining sums corresponds to generalized Reidemeister moves. As, odd writhe is an invariant under generalized Reidemeister moves and by Theorem~\ref{J(K)}, we have the upper bound for odd writhe under an arc shift move. Therefore, we have 
$$ J(L)\leq 2n.$$
Hence, $A(L)\geq J(L)/2.$ 
\end{proof}
   \begin{theorem}\label{exc}
For each positive integer $n$, there exist a 2-homogeneous proper link diagram for which arc shift is less than equal to $n$.
   \end{theorem}
   \begin{proof}
Let $D$ be the diagram of a $(2,4n)$-virtual torus link derived from the $(2,4n)$-torus link $(n\in \mathbb{N})$ by virtualizing the crossings alternately, as illustrated in Figure~\ref{rsl}. We will show $(2,4n)$- virtual torus link is a homogeneous proper link, and we will find an upper bound for its arc shift number. 
As the virtual linking numbers $L_{1}^{2}   =0$, and $L_{2}^{1}   =-2n$ are both even, by Lemma~\ref{prty}, it is arc shift equivalent to trivial virtual link.
\begin{figure}[h]
	\centering
\includegraphics[width=16cm]{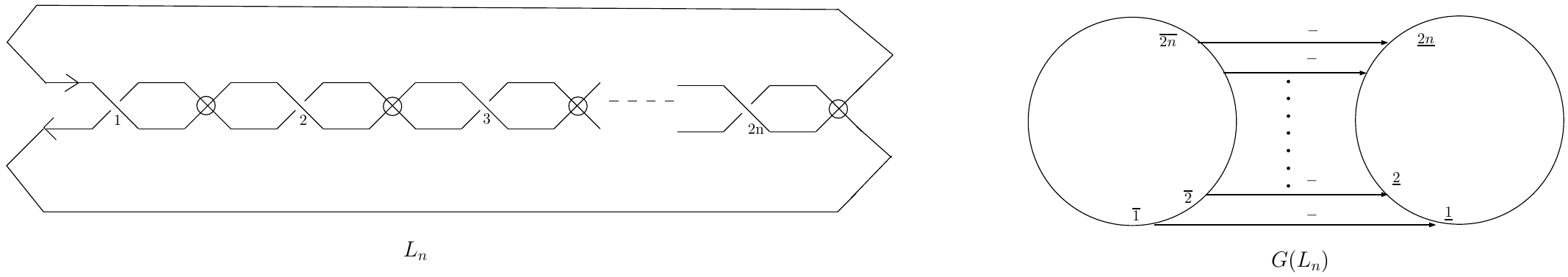}	
  \caption{ $(2,4n)$-virtual torus link diagram and the corresponding Gauss diagram.}
 \label{rsl}
   \end{figure}
For $n = 1$, it has been demonstrated that the arc shift number is 1, as shown in Figure~\ref{rst}. For $n = 2$, the $(2,8)$-virtual torus link can be trivialized by applying the arc shift move $\overline{A}_s$ first on the second real crossing and its adjacent virtual crossing, and then on the fourth real crossing and its adjacent virtual crossing, in a manner analogous to Figure~\ref{rst}. Consequently, the arc shift number $A(D) \leq 2$ for this case.
Similarly, for the $(2,4n)$-virtual torus link diagram $L$, as depicted in Figure~\ref{rsl}, we can unlink the diagram by applying $n$ arc shift moves of type $\overline{A}_s$ on the crossings $2n$ and the adjacent virtual crossings. Thus, the arc shift number of $D$ is less than or equal to $n$.

   \end{proof}

\begin{Example}
Condiser a sequence of $n$-homogeneous proper link, denoted as 
$D_n=L(K_1,\ldots,K_n)$, which is illustrated in Figure~\ref{wer}. The upper bound for the arc shift number for this virtual link diagram is $n$.
 \begin{figure}[h]
	\centering
\includegraphics[width=14cm]{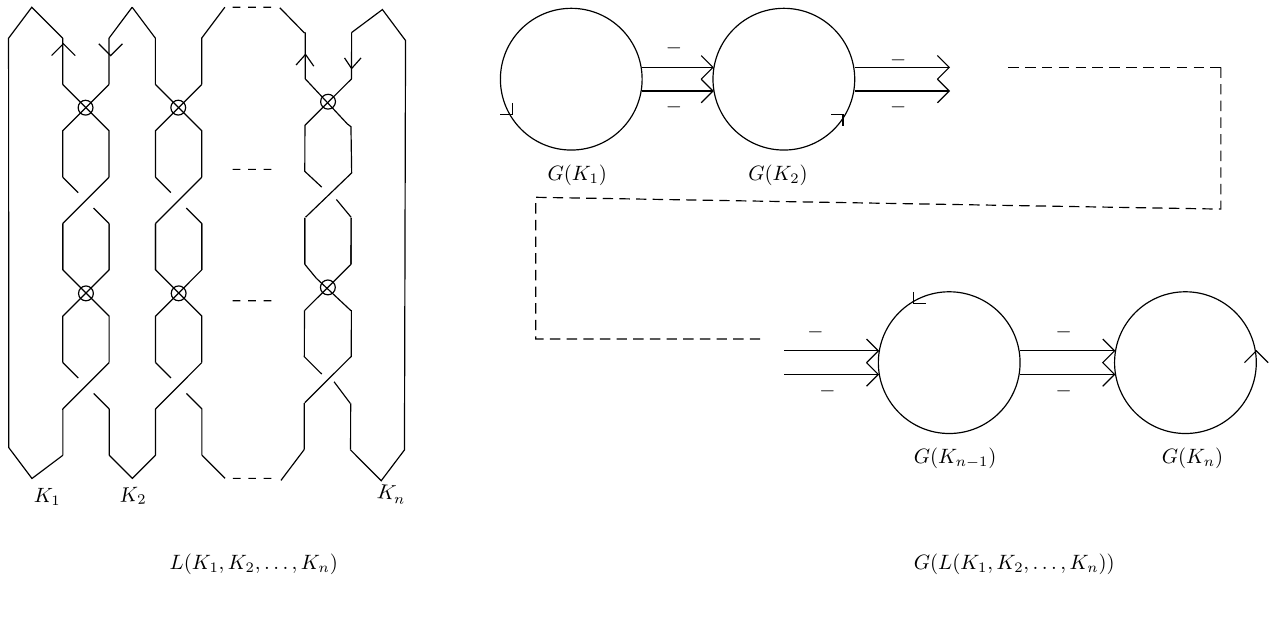}	
  \caption{The $n$-homogeneous proper link diagram $D_n$.}
 \label{wer}
   \end{figure}
\end{Example}

\begin{Example}
Consider a $2$-homogeneous proper link is shown in Figure~\ref{ctrefoil}. The calculation of the upper bound for the arc shift number of this virtual link diagram is demonstrated in Example~\ref{e1}, and the result is 3. The lower bound can be determined using Theorem~\ref{Lowerbd}. Since $J(L) = 4$, it follows that $2 \leq A(L)$. We have $2 \leq A(L) \leq 3$.
\end{Example}

\section*{Concluding remarks}

In this paper, we primarily investigate an invariant known as the arc shift number for a specific class of virtual links.
We also observe that \( n \)-component virtual links, up to arc shift equivalence, form a finite class.
Additionally, we introduce two sequences of virtual links, denoted by \( L_n \) and $D_n$, where the upper bound of the arc shift number for both \( L_n \) and $D_n$ is \( n \).
We establish that arc shift is an invariant for $n$-homogeneous proper link and we determine a lower bound for the arc shift number for this specific class of virtual links.
A compelling question for future research is whether more refined bounds can be constructed to achieve equality for the arc shift number.

\section*{Acknowledgements}
The first author expresses gratitude to IIT Ropar for the opportunity to undertake the MSc. Maths project under the supervision of Dr. M. Prabhakar.
The second author would like to thank the University Grants Commission (UGC), India, for a research fellowship with NTA Ref.No.191620008047.
The third author is deeply grateful to Professor Madeti Prabhakar for inviting him to visit and providing support and valuable encouragement about this work during his stay at IIT Ropar in May 2024.
The fourth author acknowledges the support given by the Science and Engineering Board(SERB), Department of Science $\&$ Technology, Government of India under the Mathematical Research Impact Centric Support (MATRICS) grant-in-aid with F.No.MTR/2021/00394 and by the NBHM, Government of India under grant-in-aid with F.No.02011/2/20223 NBHM(R.P.)/ R\&D II/970. 
This work was partially supported by the FIST program of the Department of Science and Technology, Government of India, Reference No. SR/FST/MS-I/2018/22(C).


\end{document}